\documentclass[11pt,letterpaper]{amsart}
\pdfoutput=1

\usepackage{amssymb}
\usepackage{amsmath}
\usepackage{amscd}
\usepackage{longtable}
\usepackage{float}
\usepackage{bigstrut}

\usepackage{graphicx}
\usepackage{multirow}
\usepackage[all,cmtip]{xy}

    \usepackage[pdftex, plainpages = false, 
		 	   pdfpagelabels,                  
		 	   pdfpagelayout = SinglePage,
		 	   bookmarks = true,
                 bookmarksopen = true,
                 bookmarksnumbered = true,
                 linktocpage,
                 pagebackref,
                 colorlinks = true,
                 linkcolor = blue,
                 urlcolor  = blue,
                 citecolor = red,
                 anchorcolor = green,
                 hyperindex = true,
                 hyperfigures]{hyperref} 

\def\@begintheorem#1#2{\par\bgroup{\sc #1\ #2. }\it \ignorespaces}
\def\@opargbegintheorem#1#2#3{\par\bgroup{\sc #1\ #2\ (#3) . }\it \ignorespaces}
\def\@endtheorem{\egroup}

\newtheorem{theorem}{Theorem}[section]
\newtheorem{corollary}[theorem]{Corollary}
\newtheorem{lemma}[theorem]{Lemma}
\newtheorem{proposition}[theorem]{Proposition}
\newtheorem{example}[theorem]{Example}
\newtheorem{definition}[theorem]{Definition}
\newtheorem{remark}{Remark}[section]
\newtheorem{notation}{Notation}[section]

\newcommand{\bt}[1]{\begin{theorem}\label{#1}}
\newcommand{\bc}[1]{\begin{corollary}\label{#1}}
\newcommand{\bl}[1]{\begin{lemma}\label{#1}}
\newcommand{\bp}[1]{\begin{proposition}\label{#1}}
\newcommand{\be}[1]{\begin{example}\label{#1}}
\newcommand{\bd}[1]{\begin{definition}\label{#1}}
\newcommand{\brem}[1]{\begin{remark}\label{#1}}
\newcommand{\bpr}{\par\noindent{\it Proof}: \ignorespaces}
\newcommand{\beq} {\begin{eqnarray}}
\newcommand{\bn}[1]{\begin{notation}\label{#1}}

\newcommand{\et}{\end{theorem}}
\newcommand{\ec}{\end{corollary}}
\newcommand{\el}{\end{lemma}}
\newcommand{\ep}{\end{proposition}}
\newcommand{\ee}{\end{example}}
\newcommand{\ed}{\end{definition}}
\newcommand{\erem}{\end{remark}}
\newcommand{\eeq}{\end{eqnarray}}
\newcommand{\en}{\end{notation}}

\def \I {{\mathcal I}}
\def \r {{\mathcal R}}
\def \q {{\mathcal Q}}

\title[Orbit closures of non-equioriented $A_3$]{\bf Resolutions of defining ideals of orbit closures for quivers of type $A_3$}

\author{\bf Kavita Sutar}

\address{Department of Mathematics\\ Northeastern University}
\email{sutar.k@husky.neu.edu}
\date{}
\keywords{Orbit closures, Cohen-Macaulay, Gorenstein, Dynkin quiver, Geometric technique, Bott's vanishing theorem}
\subjclass[2010]{14M05, 14M12, 14M17, 16G20, 16G70, 14B05, 14L30}

\begin{document}

\begin{abstract}
We investigate the properties of coordinate rings of orbit closures for quivers of type $A_3$ by considering the desingularization given by Reineke \cite{MR1985731}. We construct explicit minimal free resolutions of the defining ideals of the orbit closures thus giving us a minimal set of generators for the defining ideal. The resolution allows us to read off some geometric properties of the orbit closure. In addition, we give a characterization for the orbit closure to be Gorenstein. 
\end{abstract}
\maketitle
\section{Introduction} 
We fix an algebraically closed field $K$. Let $Q=(Q_0,Q_1)$ be a Dynkin quiver with set of vertices $Q_0$ and set of arrows $Q_1$. We use the notation $ta \stackrel{a}{\rightarrow} ha$ for arrows in $Q$.\par
A \emph{representation} $((V_i)_{i \in Q_0}, (V_a)_{a \in Q_1})$ of $Q$ is given by assigning a finite dimensional $K$-vector space $V_i$ to every vertex $i \in Q_0$ and $K$-linear maps $V_{ta} \stackrel{V_a}{\rightarrow} V_{ha}$ to every arrow $a \in Q_1$.  The \emph{dimension vector} of a representation $((V_x)_{x \in Q_0}, (V_a)_{a \in Q_1})$ is defined as the function $\underline{d}:Q_0 \longrightarrow \mathbb{Z}$ given by $\underline{d}(x)=$ dim $(V_x)$. 
Given two representations $V=((V_i)_{i \in Q_0}, (V_a)_{a \in Q_1})$ and $W=((W_i)_{i \in Q_0}, (W_a)_{a \in Q_1})$ of $Q$, a morphism $\Phi:V \rightarrow W$ is a collection of $K$-linear maps $\phi_i:V_i \rightarrow W_i$ such that for every $a \in Q_1$, the square
   \[
    \xymatrix{
    V_{ta} \ar[d]^{\phi_{ta}} \ar[r]^{V_a} & V_{ha} \ar[d]^{\phi_{ha}}\\
    W_{ta} \ar[r]^{W_a} &W_{ha}
    } 
    \]
commutes. \par
    
 With this definition of morphisms, the collection of all representations of a quiver $Q$ (over $K$) forms a category which we denote by Rep$_K(Q)$. Given a quiver $Q$, one can define its path algebra $KQ$ as the $K$-algebra generated by the paths in $Q$. It is known that $KQ$ is an associative algebra and is finite dimensional if and only if $Q$ is finite and has no oriented cycles. An important result in the theory of representation theory of associative algebras asserts that for $Q$ being a finite, connected, acyclic quiver, there is an equivalence of categories Mod $KQ$ and Rep$_K(Q)$ (refer \cite{MR2197389} for details).\par

The \emph{representation space} $Rep(Q,\underline{d})$ of a quiver $Q$ is the collection of all representations of $Q$ of fixed dimension vector $\underline{d}$.
Note that we can think of $Rep(Q,\underline{d})$ as the set $\displaystyle{\prod_{a \in Q_1} \hbox{Hom}(K^{d_{ta}}, K^{d_{ha}})}$. Thus, $Rep(Q,\underline{d})$ is a finite dimensional $K$-vector space with an affine structure.\par 

There is also the following geometric aspect to the representations of $Q$: the algebraic group $\prod_{x \in Q_0} GL(d(x))$ acts on $Rep(Q,\underline{d})$; for $V \in Rep(Q,\underline{d})$, let $\overline{O}_V$ denote the closure of an orbit $O_V$. Then $\overline{O}_V$ is a subvariety of $Rep(Q,\underline{d})$. It is an interesting problem to study the type of singularities that occur in these orbit closures. The geometry of such orbit closures was first studied by Abeasis, Del Fra and Kraft in \cite{MR626958}. They proved for the case of equioriented $A_n$ (over fields of characteristic zero) that the orbit closures are normal, Cohen-Macaulay and have rational singularities. This result was generalized to fields of arbitrary characteristic by Lakshmibai and Magyar in \cite{MR1635873}. They show using standard monomial theory that the defining ideals of orbit closures in case of equioriented $A_n$ are reduced, so the singularities of $\overline{O}_V$ are identical to those of Schubert varieties. This implies that the orbit closures are normal, Cohen-Macaulay etc.  This result was generalized to orbit closures for arbitrary quivers of type $A_n$ and $D_n$ by Bobinski and Zwara in \cite{MR1885816} and \cite{MR1967381}. They make use of certain \emph{hom-controlled} functors to reduce the general case to a special one and draw their conclusions by comparing the special case to Schubert varieties. \par

In this paper, we outline a method to contruct a minimal free resolution for the defining ideal of $\overline{O}_V$ for any Dynkin quiver $Q$. In effect, we have an algorithm for calculating the minimal resolution which depends only on the Littlewood-Richardson rule and Bott's algorithm. We refer to the method used for constructing the resolution as \emph{the geometric technique}(also referred to as the Kempf-Lascoux-Weyman geometric technique in recent literature). The general idea is to construct a desingularization $Z$ of $\overline{O}_V$ such that $Z$ is the total space of a suitable vector bundle.  Using the results of Kempf \cite{MR0384817} on collapsing of vector bundles, Lascoux \cite{MR520233} gave the construction of a minimal resolution of determinantal ideals for generic matrices. He made effective use of the combinatorics of representations of the general linear group and Bott's vanishing theorem for the cohomology of homogeneous vector bundles. These results were later generalized to similar cases. We use this generalization for our case of representations of Dynkin quivers to prove our results. A good reference for these results is the book `Cohomology of vector bundles and syzygies' by Jerzy Weyman \cite{MR1988690}.\par


In addition to giving us an explicit resolution of the coordinate ring, the geometric technique gives us a direct proof of the result of Bobinski and Zwara \cite{MR1885816} that orbit closures are normal with rational singularities in the case of non-equioriented $A_3$. The key proposition is an estimate involving the Euler form of the quiver $Q$ (Proposition \ref{prop3.1.3}).  In principle it is possible to calculate every term of the complex, although it is difficult to find a closed formula for every syzygy. However, we find a closed formula for the first term of the resolution for our case of non-equioriented $A_3$ (Theorem \ref{thm3.2.1} ). These formula allows us to calculate the minimal generators of defining ideals. We also give a characterization of Gorenstein orbits for our case (Theorem \ref{thm3.4.7}) and a sufficient condition for orbit closures to be Gorenstein for any Dynkin quiver $Q$ (Theorem \ref{thm3.4.4}). The techniques described in this paper in the context of non-equioriented $A_3$ can be generalized to other classes of Dynkin quivers. We handle these cases in our forthcoming papers.\par

 In order to find the resolution described above, we have used Reineke's desingularization \cite{MR1985731} for the orbit closure $Y$. We restrict to a 1-step desingularization in order to get semisimple vector bundles. This restriction does not induce an additional condition for non-equioriented $A_3$ since in this case every orbit admits a $1$-step Reineke desingularization.\par
This paper is organized as follows: \begin{itemize}
 \item in section $2$, we list some basic definitions and results about representations of quivers, orbit closures and
  Reineke desingularization.
 \item in section $3$, we describe the geometric setup we are working in.
 \item section $4$ contains the main results for non-equioriented $A_3$; $4.1$ contains the calculation of the 
 resolution ($F_{\bullet}$); in $4.2$ we describe the first term of $F_{\bullet}$ which gives us the minimal generators 
 of the defining ideal; in $4.3$ we investigate the last term of $F_{\bullet}$ and obtain a classification of 
 Gorenstein orbits for our case.
 \end{itemize}

\thanks{The author would like to thank Jason Ribeiro for developing the software required to calculate the Bott exchanges. It is a pleasure to thank her advisor Jerzy Weyman for suggesting the problem and for fruitful discussions. } 
\section{Preliminaries}
First, we recall some basic facts about representations of quivers.
Gabriel \cite{MR0332887} proved that the set of isomorphism classes of indecomposable representations of $Q$ is in bijective correspondence with the set of positive roots $R^+$ of the corresponding root systems. Under this correspondence, simple roots correspond to simple objects. Every representation of $Q$ can be written uniquely (upto permutation of factors) as a direct sum of indecomposable representations
  \[ V= \bigoplus_{\alpha \in R^+}m_{\alpha} X_{\alpha} \](where $m_{\alpha}=$ multiplicity of $X_{\alpha}$ in $V$). The indecomposable representations can be obtained as the vertices of the Auslander-Reiten quiver of $Q$. \par

Given a quiver $Q$, one can define an Euler form $E(Q)$ on the dimension vectors of $Q$ as follows -
\begin{definition} Let $x=(x_1,\cdots, x_n)$ and $y=(y_1,\cdots,y_n)$ be two elements of $\mathbb{N}^{Q_0}$ ($|Q_0|=n$).
Then the \emph{Euler form} $\left\langle . , . \right\rangle$ is 
  \begin{equation}
  \left\langle x , y \right\rangle=\sum_{i \in Q_0}x_iy_i - \sum_{a \in Q_1} x_{ta} y_{ha}
  \label{eulerform}
  \end{equation}
\end{definition}
 
  \brem{remark1}
  The Euler form can also be expressed in terms of the Cartan matrix $C_Q$ of $Q$ as \[\left\langle x , y   
  \right\rangle =x^t (C_Q^{-1})^t y\]
  \erem
     
  \brem{remark2} We have the following useful dimension formula in terms of the Euler form (refer \cite{MR2197389}) : if $V,W 
  \in Rep(Q,\underline{d})$, then
  \[ \left\langle \hbox{dim}~V ,~\hbox{dim}~ W \right\rangle= \hbox{dim}_K \hbox{Hom}_{KQ}(V,W)-\hbox{dim}_K \hbox{Ext}^1_{KQ}(V,W) \]
  \erem

\subsection{Orbit closures} The group $\displaystyle{\prod_{x \in Q_0}}$GL$(d(x))$ acts on Rep$(Q, \underline{d})$ by- \[((g_x)_{x \in Q_0}, (V_a)_{a \in Q_1})\longmapsto (g_{ha}V_a g_{ta}^{-1})_{a \in Q_1} . \] The orbits of this action are the isomorphism classes of representations of $Q$.

Let $V, W \in$ Rep$(Q, \underline{d})$. We say that $V \leq_{deg} W$ (i.e. $V$ is a degeneration of $W$) if the orbit of $W$ is contained in the closure of the orbit of $V$ (i.e. $O_{W} \subset \overline{O}_V$). This introduces a partial order on the orbits. There is also Riedtmann's rank criterion: $V \leq W$ if dim Hom$_Q(X,V) \leq $ dim Hom$_Q(X,W)$ for all indecomposables $X$ in $Rep(Q,\underline{d})$. The connection between these two partial orders is given by-

\bt{Bongartz' theorem}(Bongartz \cite{MR1402728}) If $A$ is a representation-directed, finite dimensional, associative $K$-algebra then the partial orders $\leq_{deg}$ and $\leq$ coincide.
\et

Since $Rep(Q,\underline{d})$ satisfies the hypotheses of this theorem, the orbit of $V \in$ $Rep(Q,\underline{d})$ is given by 
  \begin{equation}
    O_V= \{W \in Rep(Q, \underline{d}) |~ \hbox{dim Hom}_Q(X,V) = \hbox{ dim Hom}_Q(X,W) \} 
  \end{equation}
and the corresponding orbit closure is 
  \begin{equation}
    \overline{O}_V = \{W \in Rep(Q, \underline{d}) ~|~\hbox{dim Hom}_Q(X,V) \leq \hbox{ dim Hom}_Q(X,W) \}
  \end{equation}
where $X$ varies over all indecomposables in $Rep(Q,\underline{d})$.
\subsection{Desingularization}
In \cite{MR1985731}, Reineke describes an explicit method of constructing desingularizations of orbit closures of representations of $Q$. The desingularizations depend on certain directed partitions of the isomorphism classes of indecomposable objects-

\bd{Defn1} A partition $\I_*=(\I_1, \cdots , \I_s)$, where $R^+=\I_1 \cup \cdots \cup \I_s$, is called directed if:
   \begin{enumerate}
     \item $Ext^1_Q(X_{\alpha}, X_{\beta})=0$ for all $\alpha, \beta \in \I_t$ for $t=1,\cdots ,s$.
     \item $Hom_Q(X_{\beta}, X_{\alpha})=0= Ext^1_Q(X_{\alpha}, X_{\beta})$ for all $\alpha \in \I_t, \beta \in \I_u,~ t < u$
   \end{enumerate} 
\ed 

These conditions can be expressed in terms of the Euler form as-
   \begin{enumerate}
     \item $\left\langle\alpha, \beta \right\rangle=0$ for $\alpha, \beta \in \I_t$ for $t= 1, \cdots, s$
     \item $\left\langle \alpha, \beta \right\rangle \hspace{0.1cm} \geq \hspace{0.1cm} 0 \hspace{0.1cm} \geq        
            \left\langle \beta, \alpha \right\rangle$ for $\alpha \in \I_t, \beta \in \I_u, t<u$
   \end{enumerate}
 
Let $Q$ be a Dynkin quiver and let $AR(Q)$ denote its corresponding Auslander-Reiten quiver. A partition of indecomposables exists because the category of finite-dimensional representations is directed; in particular, we can choose a sectional tilting module and let $\I_t$ be its Coxeter translates.  We fix a partition $\I_*$ of $AR(Q)$. Then the indecomposable representations $X_{\alpha}$ are the vertices of $AR(Q)$. For a representation $V= \oplus_{\alpha \in R^+}m_{\alpha} X_{\alpha}$, we define representations \[ V_{(t)}:=\oplus_{\alpha \in \I_t}m_{\alpha} X_{\alpha}, ~~~~ t=1,\cdots,s \] Then $V=V_{(1)}\oplus \cdots \oplus V_{(s)}$. Let $\underline{d}_t =$ dim $V_{(t)}$. We consider the incidence variety \[Z_{\I_*, V} \subset \prod_{x \in Q_0} Flag(d_s(x),d_{s-1}(x)+d_s(x), \cdots , d_2(x)+ \cdots +d_s(x), K^{d(x)} ) \times Rep_K(Q,\underline{d}) \] defined as
  \begin{equation}
  Z_{\I_*, V} = \{((R_s(x) \subset R_{s-1}(x) \subset \cdots \subset R_2(x) \subset K^{d(x)}), V) ~|~ \forall a \in Q_1, \forall t,   
  ~~V_a(R_t(ta)) \subset R_t(ha) \}
  \label{Z}
  \end{equation}
  
In this case we say that $Z$ is a $(s-1)$-step desingularization.
\bt{Thm 1}(Reineke \cite{MR1985731}) Let $Q$ be a Dynkin quiver, $\I_*$ a directed partition of $R^+$. Then the second projection \[ q:Z_{\I_*, V} \longrightarrow Rep_K(Q, \underline{d})\] makes $Z_{\I_*, V}$ a desingularization of the orbit closure $\overline{O}_V$. More precisely, the image of $q$ equals $\overline{O}_V$ and $q$ is a proper birational isomorphism of $Z_{\I_*, V}$ and $\overline{O}_V$.
\et 

In the next section, we will realize $Z_{\I_*, V}$ as the total space of a vector bundle $\eta^*$ over \newline $\displaystyle{\prod_{x \in Q_0} Flag(d_s(x),d_{s-1}(x)+d_s(x), \cdots , d_2(x)+ \cdots +d_s(x), K^{d(x)} )}$.

\section{Geometry of orbit closures} 
We use the desingularization described above to extract useful information about orbit closures. The idea is to apply the geometric construction described in \cite{MR1988690} to calculate syzygies. We will construct a resolution $\textbf{F}_{ \bullet}$ of the coordinate ring of the orbit closure under consideration. Then using results from \cite{MR1988690} we can recover the normality of this orbit closure (we refer to \cite{MR1885816} for a different proof of the general case). In the case of $1$-step desingularization, the vector bundle $\xi$ defined below is semisimple and we get an algorithm for calculating the resolution. Also, using this resolution, we give an explicit description of the minimal set of generators for the defining ideal of the orbit closure . 
\subsection{The geometric setup} The varieties of type $Z_{\I_*, V}$ are the total spaces of homogeneous vector bundles on the product of flag varieties. We will use shorthand notation \[Z_{\underline{d}_*} \subset \prod_{x \in Q_0} Flag(\underline{d}_*{(x)}, K^{d(x)})\times Rep(Q, \underline{d}) \] to denote our incidence varieties.

The space $Rep(Q,\underline{d})$ has the structure of an affine variety. Let $A$ be the coordinate ring of $Rep(Q,\underline{d})$.
Let $\r_t(x)$ denote the tautological subbundle and $\q_t(x)$ the tautological factor bundle on $Flag(\underline{d}_*{(x)}, K^{d(x)})$. We define the following vector bundles: 
  \begin{equation}
  \xi(a) = \sum_{t=1}^s \r_t(ta)\otimes \q_t(ha)^* \subset  V(d(ta))\otimes V(d(ha))^*  
  \label{xia}
  \end{equation}
  
  \begin{equation}
  \eta(a)= V(d(ta))\otimes V(d(ha))^*/\xi(a)
  \label{etaa}
  \end{equation}
We set 
  \begin{equation}
  \eta = \bigoplus_{a \in Q_1}\eta(a)
  \label{eta}
  \end{equation}
  \begin{equation}
  \xi=\bigoplus_{a \in Q_1}\xi(a)
  \label{xi}
  \end{equation}
  
Then $Z=Z_{\underline{d}_*}$ is the total space of $\eta^*$. \vskip2mm

\noindent We have: 
    \begin{center}
    \begin{figure}[h]
    \includegraphics[scale=0.6, clip]{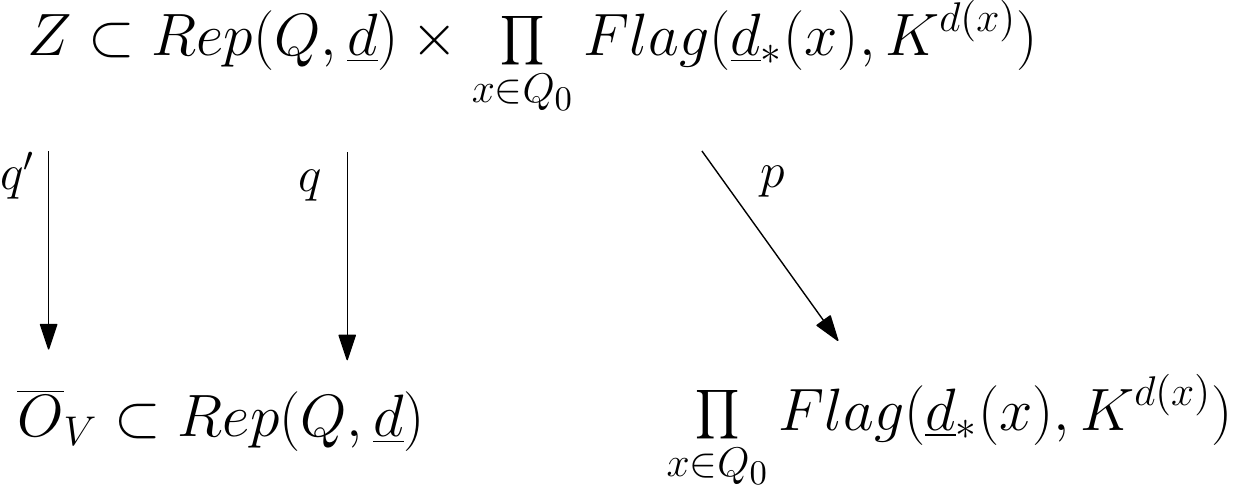}
    \label{setting}
    \end{figure}
    \end{center}

We resolve the structure sheaf $\mathcal{O}_Z$ using vector bundle $\xi$ over $\displaystyle{\prod_{x \in Q_0} Flag(\underline{d}_*{(x)},K^{d(x)})}$; this is a Koszul complex of sheaves on $Rep(Q,\underline{d}) \times \prod_{x \in Q_0}Flag(\underline{d}_*{(x)},K^{d(x)})$ -
\[ 0 \longrightarrow \bigwedge^t(p^* \xi) \longrightarrow \cdots \longrightarrow \bigwedge^2(p^* \xi) \longrightarrow p^* \xi \stackrel{\delta}{\longrightarrow} \mathcal{O}_{ Rep(Q,\underline{d}) \times \prod_{x \in Q_0}Flag(\underline{d}_*{(x)},K^{d(x)})}\]

We apply the direct image functor $R^iq_*$ to this complex to get a free resolution $\textbf{F}_{\bullet}$ of $K[\overline{O}_V]$ in terms of cohomology bundles on $V$. By Bott's theorem for general linear groups (\cite{MR1988690}, Theorem 4.1.4), the terms of this resolution are given by - \[ \textbf{F}_i = \bigoplus_{j \geq 0}H^j(\prod_{x \in Q_0}{Flag(d_*(x), V(d(x)), \bigwedge^{i+j}\xi)})\otimes A[-i-j] \]We identify $A$ with the symmetric algebra \[ \displaystyle \bigotimes_{a \in Q_1} Sym (V(ta) \otimes V(ha)^*)\]

\bt{thm2.1} \cite{MR1988690} The normalization of $\overline{O}_V$ has rational singularities if and only if $\textbf{F}_i = 0$ for $i < 0$. The orbit closure $\overline{O}_V$ is normal with rational singularities if and only if $\textbf{F}_i = 0$ for $i < 0$ and $\textbf{F}_0 = A$. 
\et

In the next section, we apply the above tool for calculations on the case of non-equioriented quiver of type $A_3$. We will consider a family of incidence varieties which is more general in the sense that we take arbitrary dimension vectors $\underline{d}_1, \cdots , \underline{d}_s$ in place of the dimension vectors described by the partition above; on the other hand we will restrict to $1$-step desingularizations. In this case, $\xi$ is semisimple and has the form 
  \begin{equation}
  \xi = \bigoplus_{a \in Q_1} \r_t(ta)\otimes \q_t(ha)^*
  \label{nicexi}
  \end{equation}

\section{Non-equioriented quiver of type $A_3$}

Let $Q$ be the non-equioriented $A_3$ quiver. We consider it in the form 
   \begin{center}
    \includegraphics[scale=0.5, clip]{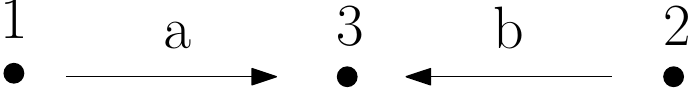}\label{dynkin}
   \end{center}
 The Auslander-Reiten quiver of $A_3$ is -
    \begin{center}
    \begin{figure}[h]
    \includegraphics[scale=0.5, clip]{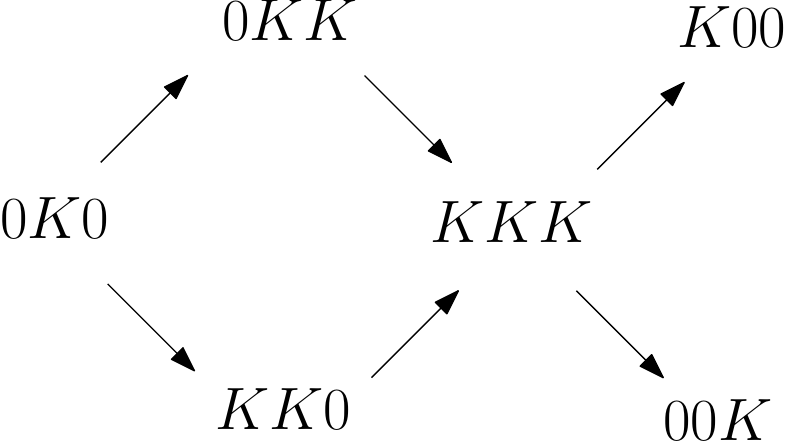}
    \caption{Auslander-Reiten quiver} \label{ARquiver}
    \end{figure}
    \end{center}
A partition of this quiver is given by choosing the sectional tilting module and its AR translates- \[ \I_1= \{K00, KKK, 00K \}, ~~\I_2= \{0K0, 0KK, KK0 \}\]
    \begin{center}
    \begin{figure}[h]
    \includegraphics[scale=0.5, clip]{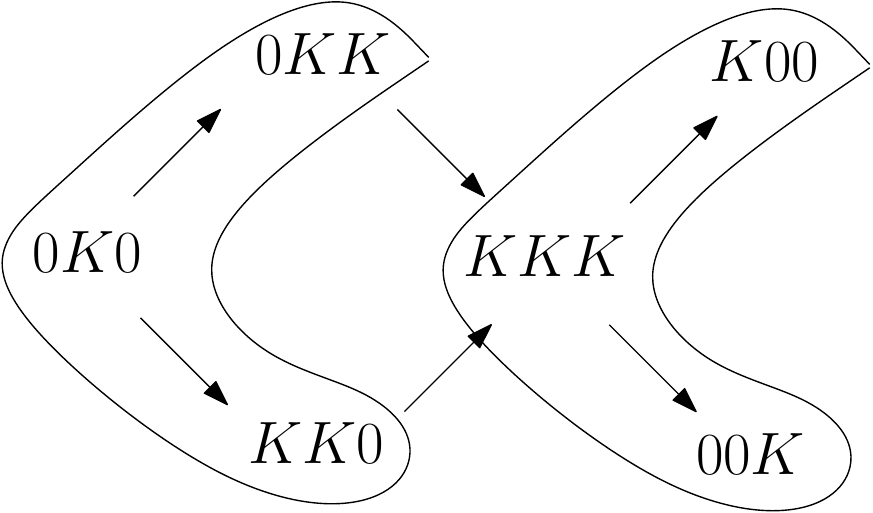}
    \caption{A partition of the AR quiver}\label{ARpart}
    \end{figure}
    \end{center}
In particular we get that there are only two parts in this partition, so every orbit admits a $1$-step desingularization. This is not true in general for other quivers and hence makes non-equioriented $A_3$ special. In general, having a $1$-step desingularization will induce additional conditions on the orbits. \par
\subsection{Calculation of $\bf F_{\bullet}$}
Let $V=a(0K0)+b(0KK)+c(KK0)+d(KKK)+e(K00)+f(00K)$ be a representation of $Q$. 
 Consider a 1-step Reineke desingularization $Z$ as follows: 
    \begin{center}
    \begin{figure}[h]
    \includegraphics[scale=0.5, clip]{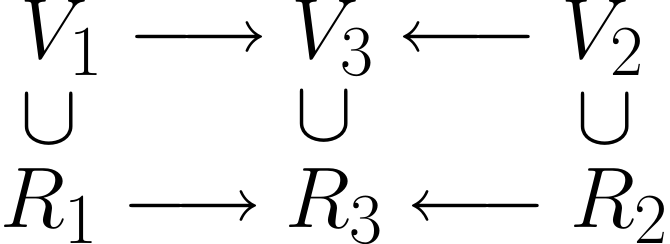}
    \caption{Desingularization of $\overline{O}_V$}
    \label{desing}
    \end{figure}
    \end{center}
 with dimension vectors of the rows being- \[\underline{d}_1+\underline{d}_2=(\alpha_1, \alpha_2, \alpha_3)\] \[\underline{d}_2=(\beta_1, \beta_2, \beta_3)\]
For example, choosing the partition in Figure \ref{ARpart} gives $\underline{d}_1+\underline{d}_2=(b+d+f,~ a+b+c+d,~ c+d+e)$ and $\underline{d}_2=(d+f,~ d,~ d+e)$.\par

Let $Q_i=V_i / R_i$, $\beta_i=$ dim $R_i$ and $\gamma_i=$ dim $Q_i$. Then the vector bundle $\xi$ is given by -
    \begin{center}
     $\xi=R_1\otimes {Q_3}^* \oplus R_2\otimes {Q_3}^*$
    \end{center}

\noindent The orbit of $V$ is \[ O_V= \{W \in Rep(Q, \underline{d}) |~ \hbox{dim Hom}_Q(X,V) = \hbox{ dim Hom}_Q(X,W),I \in \mathcal{I}_1\cup \mathcal{I}_2 \} \] and its closure is \[ \overline{O}_V= \{W \in Rep(Q, \underline{d}) |~\hbox{dim Hom}_Q(X,V) \leq \hbox{ dim Hom}_Q(X,W),I \in \mathcal{I}_1\cup \mathcal{I}_2 \} \]

\noindent The terms of the free resolution $\textbf{F}_{\bullet}$ resolving the structure sheaf are
   \begin{eqnarray}
     F_i&=&\bigoplus_{j \geq 0}H^j(\prod_{x \in Q_0}Flag(d_2(x),K^{d(x)}), \bigwedge^{i+j}\xi) \otimes A[-i-j] \\ 
        &=&\bigoplus_{j \geq 0}H^j(\prod_{x \in Q_0}Gr(d_2(x),K^{d(x)}), \bigwedge^{i+j}\xi) \otimes A[-i-j]\\ \nonumber
   \label{Fi}     
   \end{eqnarray}
\noindent Note that by Cauchy's formula we have 
   \begin{equation} 
   \bigwedge^t \xi=\bigoplus_{|\lambda|+|\mu|=t} S_{\lambda}R_1 \otimes S_{\mu}R_2 \otimes S_{\lambda'}Q_3^* \otimes  
    S_{\mu'}Q_3^*
    \label{Cauchy}
   \end{equation}

In order to calculate the cohomology terms in ($\ref{Fi}$) we apply Bott's algorithm (Theorem \ref{thm3.1.1}). The flag variety $\mathcal{V}=\prod_{x \in Q_0}Flag(d_2(x),K^{d(x)})$ is a homogeneous space for $Gl_n$. This makes it possible to describe vector bundles on $\mathcal{V}$ in terms of weights of $Gl_n$-representations (\cite{MR1988690}, Theorem 4....). We denote by $L(\alpha)$ the vector bundle corresponding to weight $\alpha$ and by $K_{\beta}$ the Weyl module corresponding to the weight $\beta$. The Bott's theorem for cohomology of vector bundles yields the following algorithm in case of $\mathcal{V}$ -

\bt{thm3.1.1} Bott's algorithm \cite{MR1988690}: Let $\alpha = (\alpha_1, . . . , \alpha_n)$. The permutation
$\sigma_i = (i, i + 1)$ acts on the set of weights in the following way:
\begin{equation} \sigma_i \cdot \alpha = (\alpha_1, . . . ,\alpha_{i-1}, \alpha_{i+1}- 1, \alpha_i + 1, \alpha_{i+2}, . . . , \alpha_n). \label{exch} \end{equation}
If $\alpha$ is a nonincreasing, then $R^0h_*L(\alpha) = K_{(\alpha)}\xi^*$ and $R^ih_*L(\alpha) = 0$ for
$i > 0$. If $\alpha$ is not a partition, then we start to apply the exchanges of type
$(\ref{exch})$, trying to move bigger number to the right past the smaller number. Two
possibilities can occur:
   \begin{enumerate}
   \item $\alpha_{i+1} = \alpha_i + 1$ when the exchange of type $(\ref{exch})$
   leads to the same sequence. In this case $R^ih_*L(\alpha) = 0$ for all $i \geq 0$.
   \item After applying say $j$ exchanges, we transform $\alpha$ into a nonincreasing
   sequence $\beta$. Then $R^ih_*L(\alpha) = 0$ for $i \neq j$ and $R^jh_*L(\alpha) = K_{(\beta)} \xi^*$
   \end{enumerate}
\et

\noindent The process of applying Bott's algorithm to weights of the form $(0^k, \alpha)$ plays an important role in all our calculations and proofs, so it is useful to introduce some notation.\par
 
\bn{notation1} Whenever we apply Bott's algorithm for the exchanges, we will refer to it as `Bott exchanges'.
\begin{enumerate}
\item We denote by $[0^k, \alpha]$ the end result after applying Bott exchanges to a weight $(0^k, \alpha)$ (here $\alpha$ is the non-increasing sequence of integers $(\alpha_1, \alpha_2, \cdots, \alpha_r)$). 
\item $N_{\alpha}$ will denote the number of Bott exchanges required to go from $(0^k, \alpha)$ to $[0^k, \alpha]$.
\end{enumerate}
\en

\noindent With this notation, applying Bott's algorithm for the weight $(0^k, \alpha)$ gives us exactly one of the following results -
\begin{enumerate}
\item[a.] during any of the successive Bott exchanges, we arrive at sequence of the form $(\cdots,m, m+1, \cdots );$ in this case, the next exchange will leave the sequence unchanged, so this is the first case of Bott's algorithm. Then we say $[0^k, \alpha]$ is zero.
\item[b.] If the above case does not occur and we reach a non-increasing sequence after $N_{\alpha}$ Bott exchanges, then we say $[0^k, \alpha] ~\hbox{is the resulting sequence}~ (\alpha_1-k, \alpha_2-k, \cdots, \alpha_{p}-k,~ p^k,~\alpha_{p+1}, \cdots, \alpha_r)$. Note that then $N_{\alpha}=pk$.
\end{enumerate}

\noindent To calculate $H^j(\mathcal{V},K^{d(x)}), \bigwedge^{i+j}\xi)$ we apply Bott's algorithm to the weights \[(0^{\gamma_1}, \lambda), (0^{\gamma_2}, \mu), (-\nu, 0^{\beta_3})\] for all representations $S_{\nu}$ occuring in $S_{\lambda'}\otimes S_{\mu'}$. 
Suppose $N_{\lambda}=uq_1$, $N_{\mu}=vq_2$ and $N_{\nu}=vr_3$. Explicitly -
 \[(0^{\gamma_1}, \lambda)=(\underbrace{0,\cdots, 0}_{\gamma_1},\lambda_1,\lambda_2,\cdots)\stackrel{Bott~ exchanges}{\longrightarrow} [0^{\gamma_1}, \lambda]=(\lambda_1-\gamma_1,\lambda_2-\gamma_1,\cdots,\lambda_u-\gamma_1,~\underbrace{u,\cdots, u,}_{\gamma_1} \lambda_{u+1},\cdots) \]
\[(0^{\gamma_2}, \mu)=(\underbrace{0,\cdots, 0}_{\gamma_2},\mu_1,\mu_2,\cdots)\stackrel{Bott~ exchanges}{\longrightarrow} [0^{\gamma_2}, \mu]=(\mu_1-\gamma_2,\mu_2-\gamma_2,\cdots,\mu_v-\gamma_2,\underbrace{v,\cdots, v,}_{\gamma_2} \mu_{v+1},\cdots)\] 
We write the third weight in its dual form- 
 \[(-\nu,0^{\beta_3})=(\cdots,-\nu_2,-\nu_1,\underbrace{0,\cdots, 0}_{\beta_3})\stackrel{Bott~ exchanges}{\longrightarrow} [-\nu,0^{\beta_3}]=(\cdots,-\nu_{w+1},\underbrace{w,\cdots, w}_{\beta_3},-\nu_w-\beta_3,\cdots,\nu_2-\beta_3,\nu_1-\beta_3)\]  Then the total number of exchanges $N$ equals $u\gamma_1+v\gamma_2+w\beta_3$. We summarize this in-

\bp{prop3.1.2} 
The terms of the complex $\textbf{F}_{\bullet}$ are given by -
\[\textbf{F}_i = \bigoplus_{t=1}^{dim \xi}\bigoplus_{|\lambda|+|\mu|=t}c_{\lambda',\mu'}^{\nu} (S_{[0^{\gamma_1},\lambda]}V_1 \otimes S_{[0^{\gamma_2},\mu]}V_2 \otimes S_{[-\nu,0^{\beta_3}]}V_3^* )\]
where $S_{\nu}\subset S_{\lambda'} \otimes S_{\mu'} $ and $|\lambda|+|\mu|- N = i$.
\ep
Since the term $|\lambda|+|\mu|- N$ occurs often, we give it a name - 
\bd{Defn} Let $\lambda$, $\mu$ and $\nu$ be partitions such that $\nu$ occurs in the Littlewood-Richardson product of 
  $\lambda$ and $\mu$. For weights $(0^p,\lambda)$,$(0^q,\mu)$ and $(-\nu,0^r)$, define  
  \[D(\lambda,\mu,\nu):=|\lambda|+|\mu|-N\] where $N=N_{\lambda}+N_{\mu}+N_{\nu}$. 
  \ed  
  \vskip0.5cm
 
From Proposition (\ref{prop3.1.2}), it is clear that in order to calculate the terms $F_i$ of the resolution, we need to calculate $D(\lambda,\mu,\nu)$. Due to the number of variables involved and the peculiar form of exchanges required, the calculation of a closed formula for $D(\lambda,\mu,\nu)$ is not easy in general. 
The next proposition is our key result which gives us a lower bound for $D(\lambda,\mu,\nu)$ in terms of the Euler form of quiver $Q$.  

\bp{prop3.1.3} $ D(\lambda,\mu,\nu) \geq \left\langle (u, v, w),(u, v, w) \right\rangle $\ep 
\bpr 
  Consider:
    \begin{eqnarray}
    \# \hbox{boxes in~} \lambda_1, \dots, \lambda_u &\geq& u^2+u\gamma_1 \nonumber \\
    \# \hbox{boxes in~} \mu_1, \dots, \mu_v &\geq& v^2+v\gamma_2 \nonumber \\
    \# \hbox{boxes in~} \nu_1, \dots, \nu_w &\geq& w^2+w\beta_3  \nonumber
    \end{eqnarray}
 Looking at the Young tableaux of $\lambda$ and $\mu$ we get (by counting boxes)-
      \begin{eqnarray}
        w.u &\geq& (\lambda_1'+ \dots + \lambda_w')-(\lambda_{u+1}+ \dots +) \nonumber \\
        w.v &\geq& (\mu_1'+ \dots + \mu_w')-(\mu_{v+1}+ \dots +) \nonumber \\ 
        \hbox{which gives~} w(u+v) &\geq& (\lambda_1'+ \dots +\lambda_w'+ \mu_1'+ \dots +\mu_w') - \nonumber\\ 
        &{}&(\lambda_{u+1}+ \dots +\mu_{v+1}+ \dots +) \nonumber \\
              &\geq& \nu_1 + \dots +\nu_w - (\lambda_{u+1}+ \dots +\mu_{v+1}+ \dots +) \nonumber \\
        \hbox{thus~} \nu_1 + \dots +\nu_w &\leq& w(u+v)+ (\lambda_{u+1}+ \dots + \mu_{v+1}+ \dots +) \nonumber          
      \end{eqnarray} 
Therefore
    \begin{eqnarray}
        (u^2+u\gamma_1)+(v^2+v\gamma_2)+(w^2+w\beta_3) &\leq& (\# \hbox{boxes in~} \lambda_1, \dots, \lambda_u) + \nonumber \\ 
        &{}& (\# \hbox{boxes in~} \mu_1, \dots, \mu_v) + \nonumber \\
        &{}& (\#\hbox{boxes in~} \nu_1, \dots, \nu_w) \nonumber \\
        &=& \lambda_1+ \dots + \lambda_u + \mu_1 + \dots +\mu_v + \nonumber\\ 
        &{}& \nu_1 + \dots + \nu_w \nonumber \\
        &\leq& \lambda_1+ \dots + \lambda_u + \mu_1 + \dots + \mu_v +  \nonumber\\
        &{}& w(u+v) + \lambda_{u+1}+ \dots + \mu_{v+1}+ \dots \nonumber\\
        &=& w(u+v) + |\lambda| + |\mu|\nonumber        
    \end{eqnarray} 
 Thus we have: 
    \begin{eqnarray}
    |\lambda| + |\mu| &\geq& (u^2+u\gamma_1)+(v^2+v\gamma_2)+w (w+\beta_3-u-v) \nonumber \\
                  &=& u\gamma_1+v\gamma_2+w\beta_3+(u^2+v^2+w^2-uw-vw) \nonumber \\
                  &=& u\gamma_1+\gamma_2+w\beta_3+\left\langle(u, v, w),(u, v, w) \right\rangle \nonumber
    \end{eqnarray} \qed    

In their paper \cite{MR1885816}, Bobinski and Zwara proved the normality of orbit closures for Dynkin quivers of type $A_n$ with arbitrary orientation. Using Proposition \ref{prop3.1.2}, we can immediately derive the normality of orbit closures in our case -

\bc{cor3.1.4} In the case of quiver $Q: 1 \rightarrow 2 \leftarrow 3$ the orbit closures are normal, Cohen-Macaulay with rational singularities.
\ec
\bpr We have that $\left\langle(u, v, w),(u, v, w) \right\rangle \geq 0$ since it the Euler form of Dynkin quiver $Q$. Then from Proposition 4.1 and 4.2, $F_i=0$ for $i < 0$. \newline \noindent Also, $\left\langle(u, v, w),(u, v, w) \right\rangle = 0$ if and only if $u=v=w=0$ in which case $\lambda=\mu=\nu=0$. Thus $F_0=0$. By Theorem 3.1, this implies that the orbit closure is normal with rational singularities.\qed

\brem{remark3} The estimate in Proposition \ref{prop3.1.2} can be extended to tree quivers having the property that every vertex is a source or sink. Thus Corollary \ref{cor3.1.4} holds for orbit closures admitting a 1-step desingularization in case of all Dynkin quivers with every vertex being source or sink. This general result is part of our next paper.
\erem

\brem{remark4} For purposes of calculation, it is useful to record some simple observations regarding the sizes of partitions $\lambda$, $\mu$ and $\nu$. From Equation (\ref{Cauchy}) it is clear that when calculating $\bigwedge^{t}\xi$, we only need to consider those partitions $\lambda$, $\mu$, $\nu$ such that $\lambda$ is contained in a $\hbox{dim}~Q_3 \times \hbox{dim}~R_1$ rectangle, $\mu$ is contained in a $\hbox{dim}~Q_3 \times \hbox{dim}~R_2$ rectangle and $\nu$ is contained in a $\hbox{dim}~(R_1+R_2) \times\hbox{dim}~Q_3 $ rectangle. Thus the largest possible contributing triples are $(\lambda, \mu, \nu)=  (\gamma_3^{\beta_1}, \gamma_3^{\beta_2}, (\beta_1+\beta_2)^{\gamma_3})$ (the notation $\alpha^\beta$ stands for the rectangular partition $(\alpha, \alpha, \dots, \alpha)$ of length $\beta$).
\erem

\be{example1}

Let $V=0K0+0KK+KK0+KKK+K00+00K$ and $I$ be the defining ideal of $\overline{O}_V$. 
Then $\alpha=(3, 4, 3)$. Let $\beta=(2, 1, 2)$. Therefore dim $\xi$ = dim $(R_1\otimes {Q_3}^* \oplus R_2\otimes {Q_3}^*)=12$. Hence we need to calculate $\wedge^{0}\xi,\wedge^{1}\xi, \dots \wedge^{12}\xi$.\par
Let $\xi_1=R_1 \otimes {Q_3}^*$ and $\xi_2=R_2 \otimes {Q_3}^*$
     \begin{eqnarray}
         \wedge^{1} \xi &=& (\wedge^{1} \xi_1 \otimes \wedge^{0} \xi_2) \oplus \wedge^{0} \xi_1 \otimes 
         \wedge^{1} \xi_2) \nonumber \\
         &=&[(S_1 R_1 \otimes S_1 {Q_3}^*)\otimes (S_0 R_2 \otimes S_0 {Q_3}^*)] \oplus [(S_0 R_1 \otimes    
         S_0 {Q_3}^*) \otimes (S_1 R_2 \otimes S_1 {Q_3}^*)] \nonumber \\
         &=& [S_1 R_1 \otimes S_0 R_2 \otimes S_1{Q_3}^* ] \oplus [S_0 R_1 \otimes S_1 R_2 \otimes 
         S_1{Q_3}^* ] \nonumber 
     \end{eqnarray}
The weight associated to the first summand is $(0,1,0; 0,0,0; 0,0,-1,0)$ and weight associated to the second summand is $(0,0,0; 0,1,0; 0,0,-1, 0)$. Applying Bott's algorithm we see that none of these terms contribute to any of the $F_i$. For an example of a contributing weight we calculate $\wedge^3 \xi$. From Remark \ref{remark4}, we know that $\lambda$ is contained in the rectangle $(3^2)$, $\mu$ is contained in $(3^2)$ and $\nu$ is contained in $(4^3)$.
     
     \begin{eqnarray} 
         \begin{split} \wedge^{3} \xi &= (\wedge^{3} \xi_1 \otimes \wedge^{0} \xi_2) \oplus (\wedge^{2} \xi_1 \otimes 
         \wedge^{1} \xi_2) \oplus (\wedge^{1} \xi_1 \otimes \wedge^{2} \xi_2) \oplus (\wedge^{0} \xi_1 \otimes \wedge^{3} 
         \xi_2) \nonumber \\
         &= [(S_{(2,1)} R_1 \otimes S_{(0)}R_2 \otimes    
           S_{(2,1)} {Q_3}^*)] \oplus [(S_{(3)} R_1 \otimes S_{(0)}R_2 \otimes S_{(1,1,1)} {Q_3}^*)]\oplus [S_{(2)}R_1 \otimes S_{(1)}  
           R_2 \otimes S_{(2,1)}{Q_3}^*)]   \\
         &\quad \oplus  [S_{(2)}R_1 \otimes S_{(1)} R_2 \otimes
           S_{(1,1,1)}{Q_3}^*)] \oplus  [S_{(1,1)}R_1 \otimes S_{(1)} R_2 \otimes S_{(2,1)}{Q_3}^*)] \oplus [S_{(1,1)}R_1 \otimes S_{(1)}
            R_2 \otimes S_{(3)}{Q_3}^*)] \\
         &\quad \oplus  [S_{(1)}R_1 \otimes S_{(2)} R_2 \otimes S_{(1,1,1)}{Q_3}^*)] \oplus [S_{(1)}R_1 \otimes S_{(2)} R_2 \otimes 
          S_{(2,1)}{Q_3}^*)] \oplus  [S_{(1)}R_1 \otimes S_{(1,1)} R_2 \otimes S_{(2,1)}{Q_3}^*)]  \\
         &\quad \oplus [S_{(1)}R_1 \otimes S_{(1,1)} R_2 \otimes S_{(3)}{Q_3}^*)] \oplus[(S_{(0)} R_1 \otimes S_{(3)}R_2 \otimes 
         S_{(1,1,1)} {Q_3}^*)]\oplus[(S_{(0)} R_1 \otimes S_{(2,1)}R_2 \otimes 
         S_{(2,1)} {Q_3}^*)]     \nonumber \\ \end{split}
     \end{eqnarray}
The weights associated to the summands in that order are:\vskip2mm
\begin{tabular}{lll}
$(0~2~1;0~0~0;0 -1 -2~0)$, & ~$(0~3~0;0~0~0;-1 -1 -1~0)$,&~$(0~2~0;0~1~0;0 -1 -2~0)$ \\
$(0~2~0;0~1~0;-1 -1 -1~0)$,&~$(0~1~1;0~1~0;0 -1 -2~0)$,&~$(0~1~1;0~1~0;0~0 -3~0)$\\
$(0~1~0;0~2~0;-1 -1 -1~0)$,&~$(0~1~0;0~2~0;0 -1 -2~0)$,&~$(0~1~0;0~1~1;0 -1 -2~0)$\\
$(0~1~0;0~1~1;0~0 -3~0)$,&~$(0~0~0;0~3~0;-1 -1 -1~0)$ ,&~ $(0~0~0;0~2~1;0 -1 -2~0)$ 
\end{tabular}\vskip2mm
Applying Bott exchanges to each weight we see that only the first and last summands contribute the non-zero terms $(\wedge^3V_1\otimes\wedge^3V_3^* \otimes A(-5))$ and $(\wedge^3V_2\otimes\wedge^3V_3^*\otimes  A(-5))$ to ${\mathbf F}_1$.
Continuing in this manner we get the resolution: 
\begin{center}
	     $\displaystyle{A=Sym(V_1 \otimes V_3^*)\oplus Sym (V_2 \otimes V_3^*)}$
	     \begin{center} $ \uparrow  $ \end{center}
	     $\displaystyle{(\wedge^3V_1\otimes\wedge^3V_3^* \otimes A(-5)) \oplus  (\wedge^3V_2\otimes\wedge^3V_3^*\otimes 
	       A(-5))\oplus (\wedge^2V_1\otimes\wedge^2V_2 \otimes \wedge^4V_3^* \otimes A(-7))} $
	   	   \begin{center} $ \uparrow  $ \end{center}
	   	 $\displaystyle{(S_{211}V_1\otimes\wedge^4V_3^*\otimes A(-6)) \oplus (S_{211}V_2\otimes\wedge^4V_3^*\otimes    
	   	 A(-6))  \oplus}$\\$\displaystyle{ (\wedge^3V_1\otimes \wedge^2V_2^* \otimes S_{2111}V_3^*\otimes A(-8)) \oplus( 
	   	 \wedge^2V_1 \otimes \wedge^3V_2 \otimes S_{2111}V_3^*\otimes A(-8))} \oplus$\\$\displaystyle{\wedge^3V_1 \otimes        \wedge^3V_2 \otimes S_{222}V_3^*\otimes A(-10)}$
	   	   \begin{center}$  \uparrow  $\end{center}  
	      $\displaystyle{(S_{211}V_1\otimes\wedge^3V_2\otimes S_{2221}V_3^*\otimes A(-11)) \oplus 
	      (\wedge^3V_1 \otimes S_{211}V_2\otimes S_{2221}V_3^*\otimes A(-11))\oplus}$\\
	      $\displaystyle{ (\wedge^2V_1\otimes S_{222}V_2^* \otimes S_{2222}V_3^*\otimes A(-13)) \oplus( 
	   	 S_{222}V_1 \otimes \wedge^2V_2 \otimes S_{2222}V_3^*\otimes A(-13))} \oplus$\\$\displaystyle{\wedge^3V_1 \otimes        \wedge^3V_2 \otimes S_{3111}V_3^*\otimes A(-9)}$
	   	   \begin{center} $ \uparrow  $ \end{center}
	     $\displaystyle{(S_{211}V_1\otimes S_{211}V_2\otimes S_{2222}V_3^*\otimes A(-12)) \oplus (S_{222}V_1 \otimes 
	     \wedge^3V_2 \otimes S_{3222}V_3^*\otimes A(-14))\oplus}$\\$\displaystyle{ (\wedge^3V_1\otimes S_{222}V_2 \otimes 
	     S_{3222}V_3^*\otimes A(-14))} $
	   	   \begin{center} $ \uparrow  $ \end{center} 
	   	 $\displaystyle{(S_{222}V_1\otimes S_{222}V_2\otimes S_{3333}V_3^* \otimes A(-19))}$
	 \end{center} 
\ee	 
\subsection{Minimal generators of the defining ideal}
Let $V \in Rep({Q, \underline{d}})$, $V=a(0K0)+b(0KK)+c(KK0)+d(KKK)+e(K00)+f(00K)$. 
Then \[\hbox{rank}~\phi = b+d,~ \hbox{rank}~\psi = c+d,~ \hbox{rank}~(\phi|\psi) = b+c+d\] We will denote these ranks by $p, q, r$ respectively. Hence $N=ub+vc+wd$. \par
We consider orbits admitting a Reineke desingularization given by the partition in Figure \ref{ARpart}. 
The following result is the main theorem of this section. It describes the first term ${\mathbf F}_1$ of the resolution ${\mathbf F}_{\bullet}$. In particular, it says that the summands of ${\mathbf F}_1$ are obtained by contributions from $\bigwedge^{rank( \phi)+1}\xi$, $\bigwedge^{rank(\psi)+1}\xi$ and $\bigwedge^{rank(\phi|\psi)+1}\xi$. As a result, we will have that the generators of the defining ideal are determinantal in the sense that they are maximal minors of $\phi$, $\psi$ and $\phi|\psi$.

\bt{thm3.2.1} ${\mathbf F}_1 = H^p(\mathcal{V},\bigwedge^{p+1}\xi) \oplus H^q(\mathcal{V},\bigwedge^{q+1} \xi) \oplus H^r(\mathcal{V},\bigwedge^{r+1} \xi)$.
\et
\bpr From Proposition \ref{prop3.1.3}, we have that 
   \[{\mathbf F}_1 = \bigoplus_{t=1}^{dim\xi}\bigoplus_{|\lambda|+|\mu|=t}c_{\lambda',\mu'}^{\nu} (S_{[0^{b},\lambda]}V_1 \otimes 
   S_{[0^{c},\mu]}V_2 \otimes S_{[-\nu,0^{d}]}V_3^* )\] 
   where $S_{\nu}\subset S_{\lambda'} \otimes S_{\mu'} $ and $D(\lambda,\mu,\nu)=1$. 
   Also by Proposition \ref{prop3.1.3}, 
      \begin{eqnarray}
      D(\lambda,\mu,\nu) &\geq& \left\langle (u, v, w),(u, v, w) \right\rangle \nonumber \\
      i.e.~~~~~~~~~~~~~~~1 &\geq& \left\langle (u, v, w),(u, v, w) \right\rangle \nonumber  
      \end{eqnarray}
   But Q is Dynkin, so $E(Q) =\left\langle (u, v, w),(u, v, w) \right\rangle > 0$, so 
   \[\left\langle (u, v, w),(u, v, w) \right\rangle =1 \]
   Thus the options for $(u,v,w)$ are $(1,0,0), (0,1,0),(0,0,1),(1,0,1),(0,1,1),(1,1,1)$. We analyze these triples to prove our     
   proposition.
   Recall that the weights of $\bigwedge^i \xi$ are of the form \[(0^b,\lambda),(0^c,\mu),(-\nu,0^d)\] where $|\lambda|+|\mu|=i$.
   
   \begin{enumerate}
   \item $(u,v,w)=(1,0,0)$: in this case $N=b$, so $|\lambda|+|\mu|= b+1$. $u=1$ implies that $\lambda=(b+1,0 \cdots, 0)$, so   
   $\mu=0$. This implies $\nu=\lambda'$, but $w=0$, so we will get a contributing triple only when $d=0$. In that case $p=\gamma_1$ and
   \[H^p(\mathcal{V},\bigwedge^{p+1}\xi) = \wedge^{p+1}V_1 \otimes \wedge^{p+1}V_3^*\] is the only contribution to ${\mathbf F}_1$.
   
   \item $(u,v,w)=(0,1,0)$: this case is analogous to the previous one. A contributing triple occurs only when $d=0$ in which case
   the contribution to ${\mathbf F}_1$ is \[H^q(\mathcal{V},\bigwedge^{q+1}\xi) = \wedge^{q+1}V_2 \otimes \wedge^{q+1}V_3^*\] 
   
   \item $(u,v,w)=(0,0,1)$: here $N=d$. So $|\lambda|+|\mu|=|\nu|=d+1$. Also $w=1$ implies $\nu$ must be $(d+1, 0,\cdots, 0)$.
   So a contributing triple occurs only when $b=c=0$. Then  $r=d$ and we get contributing triples $(1^k;~1^l;~d+1)$ where $k+l=d+1$. The 
   contribution to ${\mathbf F}_1$ is
   \[H^r(\mathcal{V},\bigwedge^{r+1}\xi) = \wedge^k V_1 \otimes \wedge^l V_2 \otimes \wedge^{r+1}V_3^*\] 
   
   \item $(u,v,w)=(1,0,1)$: this implies $N=b+d=p$. $D(\lambda,\mu,\nu)=1$ implies $|\lambda|+|\mu|-N=1$, so 
   $|\lambda|+|\mu|=|\nu|=b+d+1$. $u=1$ implies $\lambda$ is of the form $(b+1,1^k,0,\dots)$, similarly $w=1$ implies $\nu$ is of the 
   form $(d+1,1^l,0,\dots)$ (thus both $\lambda$ and $\nu$ are hooks). Then $|\nu|=b+d+1$ implies $l=b$.\par
   Since $v=0$ we know that there are zero exchanges for the weight $(0^c,\mu)$. This can 
   happen if either $\mu=0$ or $c=0$. If $\mu=0$, then $\nu=\lambda'$ and 
   \begin{eqnarray}
   H^p(\mathcal{V}, \bigwedge^{p+1}\xi)&=& S_{[0^b,\lambda]}V_1 
   \otimes S_{[-\nu,0^d]}V_3^*\nonumber \\
   &=& \wedge^{p+1}V_1 \otimes \wedge^{p+1}V_3^* \nonumber
   \end{eqnarray}
   If $c=0$, then $\mu=\nu \setminus \lambda = (1^{d-k})$.
   In this case  
   \begin{eqnarray}
   H^p(\mathcal{V}, \bigwedge^{p+1}\xi)&=& S_{[0^b,\lambda]}V_1 \otimes S_{\mu} V_2
   \otimes S_{[-\nu,0^d]}V_3^*\nonumber \\
   &=& \wedge^{b+k+1}V_1 \otimes \wedge^{d-k}V_2 \otimes \wedge^{p+1}V_3^* \nonumber
   \end{eqnarray}
 
   \item $(u,v,w)=(0,1,1)$: this case is analogous to the previous one. $u=0$ implies either $\lambda=0$ or $b=0$.
   If $\lambda=0$, then $\nu=\mu'$ and 
   \begin{eqnarray}
   H^q(\mathcal{V}, \bigwedge^{q+1}\xi)&=& S_{[0^c,\mu]}V_2 
   \otimes S_{[-\nu,0^d]}V_3^*\nonumber \\
   &=& \wedge^{q+1}V_2 \otimes \wedge^{q+1}V_3^* \nonumber
   \end{eqnarray}
   If $b=0$, then $\lambda=\nu \setminus \mu = (1^{d-k})$.
   In this case  
   \begin{eqnarray}
   H^q(\mathcal{V}, \bigwedge^{q+1}\xi)&=& S_{\lambda}V_1 \otimes S_{[0^c,\mu]} V_2
   \otimes S_{[-\nu,0^d]}V_3^*\nonumber \\
   &=& \wedge^{d-k}V_1 \otimes \wedge^{c+k+1}V_2 \otimes \wedge^{q+1}V_3^* \nonumber
   \end{eqnarray} 
   
   \item $(u,v,w)=(1,1,1)$ in this case $N=b+c+d=r$. $\lambda$ and $\mu$ are hooks of the form: \[\lambda=(b+1,1^k,0,\dots),\hskip2mm  
   \mu=(c+1,1^l,0,\dots) \]
   Since $\nu$ is such that $S_{\nu}\subset S_{\lambda'} \otimes S_{\mu'}$, $\nu$ is also a hook of the form $(d+1,1^m,0,\dots)$.
   Since $|\lambda|+|\mu|=|\nu|=b+c+d+1$, we must have $k+l=d-1$ and $m=b+c$.
   Thus
   \begin{eqnarray}
   H^r(\mathcal{V}, \bigwedge^{r+1}\xi)&=& S_{[0^b,\lambda]}V_1 \otimes S_{[0^c,\mu]} V_2
   \otimes S_{[-\nu,0^d]}V_3^*\nonumber \\
   &=& \bigoplus_{k+l=d-1}\wedge^{b+k+1}V_1 \otimes \wedge^{c+l+1}V_2 \otimes \wedge^{b+c+d+1}V_3^* \nonumber
   \end{eqnarray}
   By Cauchy's formula, this term is a direct summand of $\bigwedge^{r+1}([V_1 \oplus V_2] \otimes V_3^*)$. 
   \end{enumerate}
\qed

\bc{cor3.2.2} Let rank $(\phi)=p$, rank $(\psi)=q$, rank $(\phi+\psi)=r$. The minimal generators of the defining ideal are determinantal: $(p+1) \times (p+1)$ minors of $\phi$, the $(q+1) \times (q+1)$ minors of $\psi$ and the $(r+1) \times (r+1)$ minors of $\phi|\psi$, taken by choosing $b+k+1$ columns of $\phi$ and $c+l+1$ columns of $\psi$, where $k+l=d-1$.
\ec
\bpr The defining ideal of the orbit closure $\overline{O}_V$ is generated by the image of the map $F_1 \stackrel{\delta}{\longrightarrow} A$. By Theorem \ref{thm3.2.1}, the image of the differential map $\delta$ is generated by $(p+1)\times (p+1)$-minors of the matrix corresponding to the linear map $\phi$, $(q+1)\times (q+1)$-minors of the matrix corresponding to the linear map $\psi$ and $(r+1)\times (r+1)$-minors of the matrix corresponding to the linear map $\phi|\psi$. \qed
\vskip0.5cm
 
 In Example \ref{example1}, we found
	   \[{\mathbf F}_1 = (\wedge^3V_1\otimes\wedge^3V_3^* \otimes A(-5)) \oplus  (\wedge^3V_2\otimes\wedge^3V_3^*\otimes 
	       A(-5))\oplus (\wedge^2V_1\otimes\wedge^2V_2 \otimes \wedge^4V_3^* \otimes A(-7))\]
Fixing a basis for vector spaces $V_1$, $V_2$ and $V_3$, the minimal generators of the defining ideal are $3 \times 3$ minors of the $4 \times 3$ matrices $\phi$ and $\psi$ and $4 \times 4$	minors of the map $\phi | \psi : V_1 \oplus V_2 \rightarrow V_3$, obtained by choosing $2$ columns of $\phi$ and $2$ columns of $\psi$.

\subsection{ ${\mathbf F}_{top}$ and classification of Gorenstein orbits}

Lets consider a Dynkin quiver $Q$. We denote the last term of the resolution ${\mathbf F}_{\bullet}$ by ${\mathbf F}_{top}$. Let $t=\hbox{dim}~\xi$, where $\xi$ is the vector bundle defined in (\ref{xi}). The top exterior power of $\xi(a)$ is 
\begin{eqnarray}  S_{[0^{d_1(ta)}, d_1(ha)^{d_2(ta)}, \cdots, (d_1(ha)+\cdots+d_{s-1}(ha))^{d_s(ta)}]}(ta) \\ \nonumber \otimes S_{[(-d_2(ta)-\cdots-d_s(ta))^{d_1(ha)}, \cdots, -d_s(ta)^{d_{s-1}(ha)}, 0^{d_s(ha)}]}(ha)^*\end{eqnarray}
Thus the top exterior power of $\xi$ is given by 
   \begin{equation} \bigwedge^{t} \xi= \bigotimes_{x \in Q_0} S(x)(k_1(x)^{d_1(x)},\cdots,k_s(x)^{d_s(x)} )   
   \label{topterm}
   \end{equation}
   where 
   \begin{equation} k_p(x)= \sum_{a \in Q_1; ta=x} \sum_{u<p} d_u(ha)-\sum_{a \in Q_1; ha=x} \sum_{u>p} d_u(ta)
   \end{equation} \par

First, we give a sufficient condition for the orbit closure $\overline{O}_V$ to be Gorenstein in case of any Dynkin quiver $Q$. 
The condition that for every $x \in Q_0$, the number 
   \begin{equation} k_p(x)-\sum_{u<p}d_u(x)+\sum_{u>p}d_u(x) \label{condition}
   \end{equation} 
is independent of $p$ $(p=1,2, \cdots, s)$, is equivalent to the the condition that $\bigwedge^t \xi$, the top exterior power of $\xi$, contributes a trivial representation to ${\mathbf F}_{top}$. We show that the latter condition, together with normality, implies that the corresponding orbit closure is Gorenstein.
First we show that the condition (\ref{condition}) is equivalent to the property that the $\tau$-orbits in the Auslander-Reiten quiver are constant- 

\bl{lemma3.4.1} Let $\tau$ denote the Auslander-Reiten translate and suppose $\underline{d}(x)=(d_u(x))_{u=1,2,\cdots, s}$ are dimensions of the flag at vertex $x$ in the desingularization $Z$. Then
\[ \langle \underline{e}_x, \underline{d}_p(x) \rangle = -\langle \underline{d}_{p+1}(x), \underline{e}_x \rangle  \] 
for all $x \in Q_0$ and $p=1,2,\cdots, s-1$, where $\underline{e}_x$ is the dimension vector of the simple representation supported at $x$.
\el 
   \bpr (\ref{condition}) translates to the equations- 
   \begin{equation} k_{p+1}(x)-k_p(x)=d_p(x)+d_{p+1}(x)
   \end{equation}
   for $x \in Q_0$ and $p=1,2,\cdots, s-1$. This is equivalent to-   
   \begin{equation} \sum_{a \in Q_1;ta=x}d_t(ha)+\sum_{a \in Q_1; ha=x}d_{p+1}(ta)=d_{p+1}(x)+d_p(x)
   \end{equation}
   for all $x \in Q_0$ and $p=1,2,\cdots, s-1$. These conditions can be expressed in terms of Euler form as follows-
   \begin{eqnarray} \nonumber
   \langle \underline{e}_x, \underline{d}_p \rangle &=& d_p(x)-\sum_{\substack{a \in Q_1 \\ ta=x}}d_p(ha) \\ \nonumber
   &=& \sum_{\substack{ a\in Q_1 \\ ha=x}}d_{p+1}(ta)-d_{p+1}(x) \\ \nonumber
   &=& -\langle \underline{d}_{p+1}, \underline{e}_x \rangle
   \end{eqnarray} 
   Thus, 
   \begin{equation} \langle \underline{e}_x, \underline{d}_p \rangle = -\langle \underline{d}_{p+1}, \underline{e}_x
    \rangle \label{euler}
   \end{equation}
   where $\underline{e}_x$ is the dimension vector of the simple representation supported at $x$.   
   \qed

\bl{lemma3.4.2} Let $m=$ dim $\mathcal{V}$ and $t =$ dim $\xi$. Then \[\hbox{codim}~ \overline{O}_V = t-m\]
\el
   \bpr
   \begin{eqnarray}
   \hbox{codim}~\overline{O}_V &=& \hbox{dim}~ X - \hbox{dim}~ \overline{O}_V  \nonumber\\
   &=& \hbox{dim}~ X - \hbox{dim}~ Z  \nonumber\\
   &=& \hbox{dim}~ X -( \hbox{dim}~ X+m-t)  \nonumber\\
   &=& t-m \nonumber
   \end{eqnarray}
   \qed

\bl{lemma3.4.3} Suppose $\bigwedge^t \xi$ contributes a trivial representation to ${\mathbf F}_{t-m}$. Then the resolution ${\mathbf F}_{\bullet}$ is self-dual. In particular, ${\mathbf F}_{t-m} \cong {\mathbf F}_0^*$.
\el
   \bpr  If $H^m(\mathcal{V},\bigwedge^t \xi)$ is a trivial representation then $\bigwedge^t \xi \cong \omega_{\mathcal{V}}$, where 
   $\omega_{\mathcal{V}}$    
   denotes the canonical sheaf on $\mathcal{V}$. This implies that $\omega_{\mathcal{V}} \otimes \bigwedge^t \xi^*  \cong \bigwedge^0 \xi 
   \cong K$. Then for $0 \leq i \leq m$,
   \begin{eqnarray}
    {\mathbf F}_{t-m-i}&=& \bigoplus_{j \geq 0} H^{m-j}(\mathcal{V}, \bigwedge^{t-i-j}\xi) \\ \nonumber
    &\cong& \bigoplus_{j \geq 0} H^j(\mathcal{V}, \omega_{\mathcal{V}} \otimes \bigwedge^{t-i-j}\xi^*)^* ~~~~~~~~~~(\hbox{by Serre 
    duality})\\  
    \nonumber   
    &\cong& \bigoplus_{j \geq 0} H^j(\mathcal{V}, \omega_{\mathcal{V}} \otimes \bigwedge^t \xi^* \otimes \bigwedge^{i+j} \xi)^*\\ 
    \nonumber
    &\cong& \bigoplus_{j \geq 0} H^j(\mathcal{V}, \bigwedge^{i+j} \xi)^* \\ \nonumber
    &=&  {\mathbf F}_i^* 
   \end{eqnarray}
   \qed

\bt{thm3.4.4} Assume that for each $p=1,2, \cdots, s-1$ we have $\underline{d}_{p+1}=\tau^+ \underline{d}_p$. Then the complex ${\mathbf F}_{\bullet}$ is self-dual. If the incidence variety comes from Reineke desingularization and the corresponding orbit closure is normal with rational singularities, then it is also Gorenstein. 
\et
  \bpr If the $\tau$-orbits of an AR quiver are constant then by Lemma \ref{lemma3.4.1}, $\bigwedge^t \xi$ contributes a trivial 
  representation to ${\mathbf F}_{t-m}$. Then applying Lemma \ref{lemma3.4.3} we get that ${\mathbf F}_{t-m} \cong {\mathbf F}_0^* \cong 
  A^*$, therefore dim ${\mathbf F}_{t-m}=1$.
  \qed

In particular, for our case of non-equioriented $A_3$ this says that the orbits with multiplicities satisfying $a=d$, $b=e$ and $c=f$
are Gorenstein. \par

Next, we investigate necessary conditions for for the orbit closure $\overline{O}_V$ to be Gorenstein in case of non-equioriented $A_3$.
Recall that for our case of non-equioriented $A_3$, we have desingularization-
  \begin{center}
  \includegraphics[scale=0.5, clip]{sec3/desing.pdf}
  \end{center}
As before, let $V=a(0K0)+b(0KK)+c(KK0)+d(KKK)+e(K00)+f(00K)$ be a representation of $A_3$.   
Then \[d_1=(b, ~a+b+c,~ c) ;~~ d_2=(d+f, ~d,~ d+e)\]  
From (\ref{topterm}) the weights for $\bigwedge^t \xi$ are:
\[(0^b, (a+b+c)^{d+f}),~~(0^c, (a+b+c)^{d+e}),~((-2d-e-f)^{a+b+c},0^d) \]

For the case of non-equioriented $A_3$, we investigate the following question: in what cases does $\bigwedge^{t}\xi$ contribute a non-zero representation? To which term $F_i$ does $\bigwedge^{t}\xi$ contribute?  

First we show that a contribution from $\bigwedge^{t}\xi$ always goes to ${\mathbf F}_{t-m}$.
\bl{lemma3.4.5} If the weight of the $\bigwedge^{t}\xi$ gives a non-zero partition after Bott exchanges, then the corresponding representation is a summand of ${\mathbf F}_{t-m}$. 
\el
   \bpr It is enough to show that $D(\lambda,\mu,\nu)=$ codim $\overline{O}_V$ for $\lambda=( (a+b+c)^{d+f})$ and $\mu=((a+b+c)^{d+e})$.
   We apply Bott's algorithm to each weight to get:
   \[[0^b, (a+b+c)^{d+f}] = ((a+c)^{d+f}, (d+f)^b) ~\hbox{after}~  b(d+f) ~\hbox{Bott exchanges},\]
   \[[0^c, (a+b+c)^{d+e}] = ((a+b)^{d+e}, (d+e)^c) ~\hbox{after}~ c(d+e) ~\hbox{Bott exchanges},\]
   \[[(-2d-e-f)^{a+b+c},0^d] = ((-a-b-c)^d, (-d-e-f)^{a+b+c}) ~\hbox{after}~ d(a+b+c) ~\hbox{Bott exchanges}.\] 
   \begin{eqnarray} \nonumber
   D(\lambda, \mu)&=&[(d+f)(a+b+c)]+[(d+e)(a+b+c)]\\ \nonumber
   &&-[b(d+f)+c(d+e)+d(a+b+c)] \\ \nonumber
   &=&ad+ae+af+be+cf \\ \nonumber
   &=&\hbox{codim}~ \overline{O}_V \\ \nonumber
   &=&t-m
   \end{eqnarray}    
   \qed

Next we list the cases in which $\bigwedge^{t}\xi$ contributes a non-zero term. Observe that a contribution will occur whenever the Bott exchanges give a non-increasing sequence for every term of
\[ (0^b, (a+b+c)^{d+f}),~~(0^c, (a+b+c)^{d+e}),~((-2d-e-f)^{a+b+c},0^d)  \]
Also, note that if any of $b,c$ or $d$ are zero, then there are no exchanges for the corresponding term in the weight. We base our cases on this observation.

\bp{Prop3.4.6} $\bigwedge^{t}\xi$ contributes to ${\mathbf F}_{t-m}$ in the following cases when the corresponding conditions are satisfied:
    \begin{table}[h]
    \centering
    \renewcommand{\arraystretch}{1.3}
    \begin{tabular}{|c|ccc|}
    \hline
    Cases & \multicolumn{3}{|c|}{Conditions}\\
    \hline
    $b=0, ~c=0, ~d=0$ &   no condition &&      \\
    \hline
    $b\neq 0, ~c=0, ~d=0$ &  $a+c \geq d+f$ &  &      \\
    \hline
    $b=0, ~c\neq0, ~d=0$ &  & $a+b \geq d+e$ &  \\
    \hline
    $b=0,~ c=0,~ d\neq0$ &  &  &$ d+e+f \geq a+b+c  $      \\
    \hline
    $b=0, ~c\neq0,~ d\neq0$ &  & $a+b \geq d+e$ ,&~ $d+e+f \geq a+b+c$  \\
    \hline
    $b\neq 0, ~c=0, ~d\neq0$ & $a+c \geq d+f$ ,& & $ d+e+f \geq a+b+c $       \\
    \hline
    $b\neq 0,~ c\neq0, ~d=0$ &$ a+c \geq d+f$, &~$a+b \geq d+e$ & \\
    \hline
    $b\neq0,~c\neq0,~d\neq0$ & $a+c \geq d+f $, & ~$a+b \geq d+e $, &~$ d+e+f \geq a+b+c$ \\
    \hline
    \end{tabular}
    \caption{}
    \label{table1}
    \end{table}
\ep

For the cases listed above, we calculate the representation that $\bigwedge^t \xi$ contributes to ${\mathbf F}_{t-m}$: 
    \begin{center}
    \begin{longtable}{|c|c|c|}
    \hline
    Case & Weight of $\bigwedge^t \xi$ & Corresponding term $H^m(\mathcal{V},\bigwedge^t \xi )$ in ${\mathbf F}_{t-m}$ \bigstrut\\
    \hline
       \endfirsthead
       \hline 
          \multicolumn{1}{|c|}{Case} &          
          \multicolumn{1}{|c|}{Weight of $\bigwedge^t \xi$} &
          \multicolumn{1}{|c|}{Corresponding term $H^m(\mathcal{V},\bigwedge^t \xi)$ in ${\mathbf F}_{t-m}$} \bigstrut \\
      \hline 
      \endhead
      \hline 
           \multicolumn{3}{|r|}{{Continued on next page}} \\ \hline
      \endfoot
      \endlastfoot
    $b=0,~c=0,~d=0$ & $ (a^{f};~~a^{e};~(-e-f)^{a}) $&  $S_{(a^{f})}V_1 \otimes S_{(a^{e})}V_2 \otimes S_{((-e-f)^{a})}V_3^*$     
    \bigstrut\\
    \hline
    $b\neq 0,~c=0,~d=0$ & $ (0^b,(a+b)^{f};~~(a+b)^{e};~(-e-f)^{a+b}) $&  $S_{(a^{f},f^b)}V_1 \otimes S_{((a+b)^{e})}V_2 \otimes   
    S_{((-e-f)^{a+b})}V_3^*$ \bigstrut\\
    \hline
    $b=0, ~c\neq0, ~d=0$ &$ ((a+c)^{f};~~0^c,(a+c)^{e};~(-e-f)^{a+c}) $&  $S_{((a+c)^{f})}V_1 \otimes S_{(a^{e}, e^c)}V_2 \otimes  
    S_{((-e-f)^{a+c})}V_3^*$  \bigstrut \\
    \hline
    $b=0,~ c=0,~ d\neq0$ & $ (a^{d+f};~~a^{d+e};$&  $S_{(a^{d+f})}V_1 \otimes S_{(a^{d+e})}V_2$  \\
    &$~(-2d-e-f)^{a},0^d) $ & $ \otimes S_{(-a^{d},(-d-e-f)^{a})}V_3^*$ \bigstrut\\
    \hline
    $b=0, ~c\neq0,~d\neq0$ &  $ ((a+c)^{d+f};~~0^c,(a+c)^{d+e};$&  $S_{((a+c)^{d+f})}V_1 \otimes S_{(a^{d+e},(d+e)^c)}V_2 $   \\
    & $~(-2d-e-f)^{a+c},0^d) $& $\otimes S_{((-a-c)^{d},(-d-e-f)^{a+c})}V_3^*$ \\
    \hline
    $b\neq 0, ~c=0, ~d\neq0$ &  $ (0^b,(a+b)^{d+f};~~(a+b)^{d+e};$&  $S_{(a^{d+f},(d+f)^b)}V_1 \otimes S_{((a+b)^{d+e})}V_2 $  \\
    &$~(-2d-e-f)^{a+b},0^d) $& $\otimes S_{((-a-b)^{d},(-d-e-f)^{a+b})}V_3^*$ \\
    \hline
    $b\neq 0,~ c\neq0, ~d=0$ &  $ (0^b,(a+b+c)^{f};~~0^c,(a+b+c)^{e};$ &  $S_{((a+c)^{f},f^b)}V_1 \otimes S_{((a+b)^{e},e^c)}V_2   $ \\
    &$~(-e-f)^{a+b+c}) $ & $\otimes S_{((-e-f)^{a+b+c})}V_3^*$ \\
    \hline
    $b\neq0,~c\neq0,~d\neq0$ &  $ (0^b,(a+b+c)^{d+f};~~0^c,(a+b+c)^{d+e};$&  $S_{((a+c)^{d+f},(d+f)^b)}V_1 \otimes  
    S_{((a+b)^{d+e},(d+e)^c)}V_2$   \\
     & $~(-2d-e-f)^{a+b+c}, 0^d) $ &  $\otimes S_{((-a-b-c)^{d},(-d-e-f)^{a+b+c})}V_3^*$\\
    \hline
    \caption{Contribution from $\bigwedge^t \xi$} \label{table2}
    \end{longtable}
    \end{center}
 
Since $\overline{O}_V$ is Cohen-Macaulay by Corollary \ref{cor3.1.4}, it is Gorenstein if and only if ${\mathbf F}_{t-m}$ is $1$-dimensional. It is known from the work of Lascoux that determinantal varieties are Gorenstein. So if an orbit closure can be viewed as a determinantal variety, it will be Gorenstein. We list such cases after Theorem \ref{thm3.4.7}.
 
\bt{thm3.4.7} A non-determinantal orbit closure $\overline{O}_V$ is Gorenstein if and only if $V$ is in an orbit with multiplicities satisfying one of the following properties:
  \begin{enumerate}
  \item[(1)] $a=d$, $b=e$, $c=f$
  \item[(2)] $a=d+e$, $b=0$, $c=f$
  \item[(3)] $a=d+e$, $b=f=0$
  \item[(4)] $a=d+f$, $c=0$, $b=e$
  \item[(5)] $a=d+f$, $c=e=0$ 
  \end{enumerate}
\et
   \bpr Part $(1)$ follows from Theorem \ref{thm3.4.4} and Table \ref{table2}. For instance, in the case $b\neq0,~c\neq0,~d\neq0$ 
   the term $H^m(\mathcal{V},\bigwedge^t \xi)$ is 1-dimensional if and only if $a+c=d+f$, $a+b=d+e$ and $a+b+c=d+e+f$ that is if and only 
   if $a=d$, $b=e$ and      $c=f$.  For the remaining parts, note that $(2)$ is symmetric to $(4)$ and $(3)$ is symmetric to $(5)$, so it 
   suffices to prove $(2)$ and $(3)$.\par
   For part $(2)$, note that the weight of $\bigwedge^t \xi$ is \[((d+e+c)^{d+c};~0^c,(d+e+c)^{d+e};~(-2d-e-c)^{d+e+c},0^d)\] Calculating 
   $D(\lambda,\mu, \nu)$ shows that $H^m(\mathcal{V}, 
   \bigwedge^t \xi)$ is non-zero and dim $H^m(\mathcal{V}, \bigwedge^t \xi)=1 $. So by Lemma \ref{lemma3.4.3}, the complex ${\mathbf 
   F}_{\bullet}$ 
   is self-dual in this case. ${\mathbf F}_0=A$ implies ${\mathbf F}_{t-m}$ is $1$-dimensional, hence Gorenstein. \par
   Finally, to prove part $(3)$ we show combinatorially that there exists a unique triple $(\lambda,\mu,\nu)$ for which 
   $D(\lambda,\mu,\nu)=t-m$. Notice that for this case we have $t-m = (d+e+c)(2d+e)- d(d+e+c)-c(d+e)= (d+e)^2$.\par
   \begin{itemize}
   \item[Claim 1:] $D((d+e)^d;(d+e+c)^{d+e};(2d+e)^{d+e},(d+e)^c)=t-m$. \par
      \begin{eqnarray}\nonumber
      D((d+e)^d;(d+e+c)^{d+e};(2d+e)^{d+e},(d+e)^c) &=& (d+e)(2d+e+c)-c(d+e)-d(d+e)\nonumber \\ 
      &=& (d+e)^2 \nonumber\\ \nonumber
      \end{eqnarray} 
      Also note that $(2d+e)^{d+e},(d+e)^c$ is the unique term in the Littlewood-Richardson product of $(d+e)^d$ and $(d+e+c)^{d+e}$ 
      which  satisfies conditions ...   
   \item[Claim 2:] If $(\hat{\lambda}, \hat{\mu}, \hat{\nu})$ is any other contributing triple, then $D(\hat{\lambda}, \hat{\mu}, 
   \hat{\nu}) < t-m$. \par
   Observe that $\nu$ has $2$ corner boxes either of which can be removed to obtain a smaller $\hat\nu$. Suppose we remove the first 
   corner box. This corresponds to removing one corner box from $\mu$. The next triple contributing a $1$-dimensional representation is 
   $(\hat{\lambda}, \hat{\mu}, \hat{\nu}) = ((d+e-1)^d;~(d+e+c)^{d+e-1}, d+e-1;~(2d+e-1)^{d+e-1},(d+e-1)^{c+1})$ 
   with number of exchanges decreased by $c+d$. Then
      \begin{eqnarray}\nonumber
      D(\hat{\lambda}, \hat{\mu}, \hat{\nu}) &=& (d+e-1)(2d+e+c-1)- c(d+e-1)+d(d+e-1) \nonumber\\ 
      &=& (d+e-1)^2 < t-m \nonumber\\ \nonumber
      \end{eqnarray}
   On the other hand if we remove the second corner box, this corresponds to removing a box from $\mu$ and the next contributing triple 
   is again $((d+e-1)^d;~(d+e+c)^{d+e-1}, d+e-1;~(2d+e)^{d+e-1},(d+e-1)^{c+1})$. Thus, removing boxes from either corner results in a  
   triple with $D(\hat{\lambda}, \hat{\mu}, \hat{\nu}) < t-m$.    
   \end{itemize}
   Thus, the $(d+e)^d;(d+e+c)^{d+e};(2d+e)^{d+e},(d+e)^c)$ is the unique triple that contributes to $F_{t-m}$; applying Bott's exchanges 
   to the corresponding weight we get that the contribution is a trivial representation. By Lemma \ref{lemma3.4.3} and the fact that 
   $\overline{O}_V$ is normal, we're done.
   \qed 
   \vskip0.5cm    
   
Finally, we give a list of orbits that can occur if the orbit closure is Gorenstein. These are the determinantal orbits mentioned earlier. Since it is enough to specify the multiplicities $a, b, c, d, e, f$ to specify an orbit, we present the orbits in the shape of the AR quiver (Figure \ref{ARquiver}) with multiplicities in place of indecomposables.
\begin{center}
  \begin{figure}[H]
  \includegraphics[scale=0.6, clip]{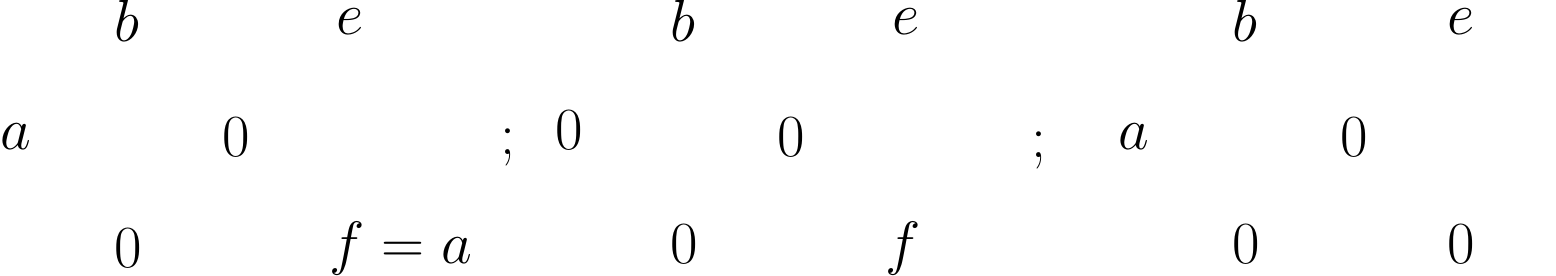}
  \label{orbits1}
  \end{figure}
  
  \begin{figure}[H]
  \includegraphics[scale=0.6, clip]{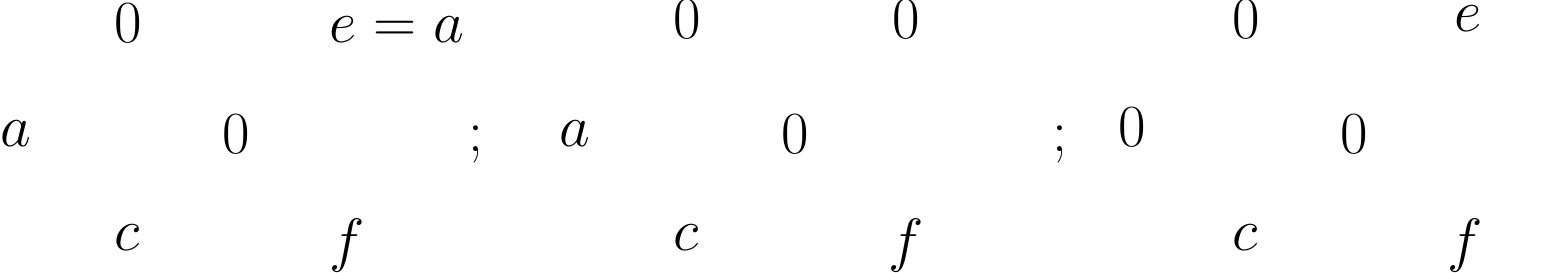}
  \label{orbits2}
  \end{figure}

  \begin{figure}[H]
  \includegraphics[scale=0.6, clip]{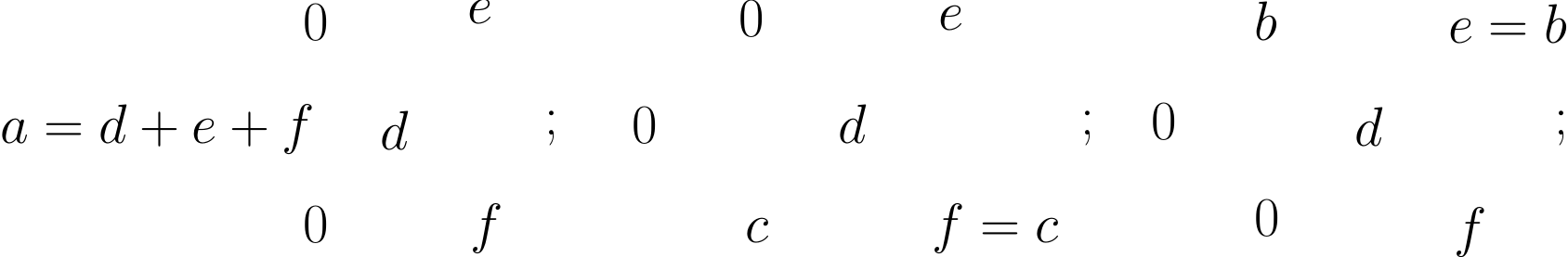}
  \label{orbits3}
  \end{figure}

  \begin{figure}[H]
  \includegraphics[scale=0.6, clip]{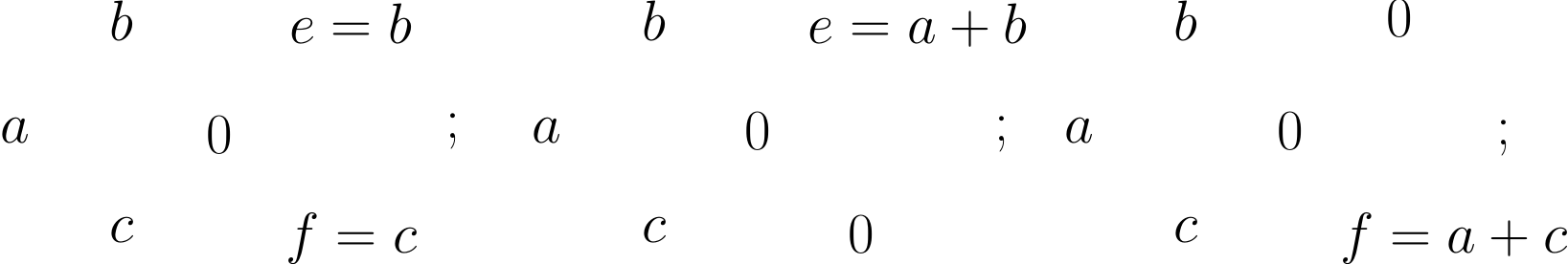}
  \label{orbits4}
  \end{figure}

  \begin{figure}[H]
  \includegraphics[scale=0.6, clip]{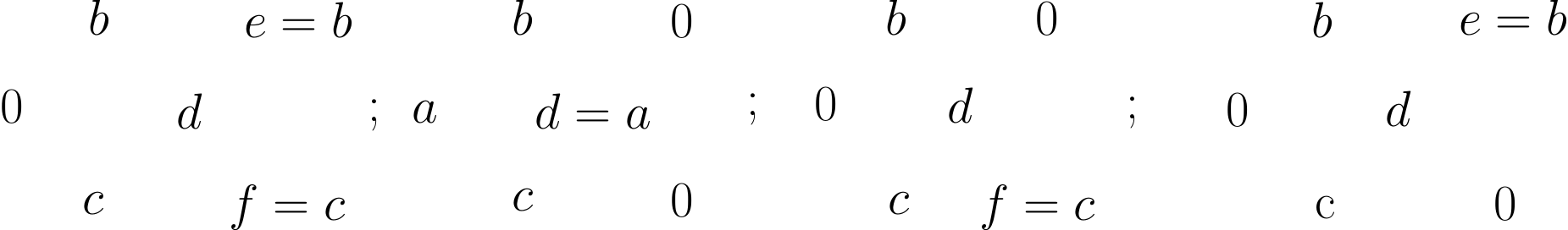}
  \label{orbits5}
  \end{figure}
\end{center}  
   

\bibliographystyle{alpha}
\bibliography{references3}

\end{document}